Analysis.

## About the operator creating secondary polynomials.

### Roland Groux

*Lycée Polyvalent Rouvière, rue Sainte Claire Deville,
BP 1205. 83070 Toulon .Cedex. France.*

Email : roland.groux@orange.fr

**Abstract.**

We are studying here the classical operator creating secondary polynomials associated with an orthogonal system for a continuous probability density function on a real interval. We know it is possible with the coupling of Stietjes Transforms to build an auxiliary measure also called secondary measure, making associated polynomials orthogonal. One of the consequences of this definition is the possibility to extend this operator to a function having interesting isometric characters. [3].

Under some hypotheses, there appears in the construction of the secondary measure a function having a privileged role which we will develop here: we will call it the reducer which will enable us among other kings to reverse the studied extension and to make its adjunct explicit. We will also illustrate it, in the case of classical orthogonal polynomials with a few results on the Fourier's coefficients of this reducer with different numerical applications.

______________________________________________________________________

## 1. Introduction and notations.

In what follows we consider a probability density function $x \mapsto \rho(x)$ on an interval $I$ bounded by $a$ and $b$. $L^2(I,\rho)$ is the associated Hilbert'space, provided with the classical inner product : $(f,g) \mapsto <f/g>_\rho = \int_a^b f(t)g(t)\rho(t)dt$ and Hilbert's basis of normalized polynomials written as : $n \mapsto P_n$. Let us call $Q_n(X) = \int_a^b \frac{P_n(t) - P_n(X)}{t - X}\rho(t)dt$ the so-called secondary polynomial associated with $P_n$.

We also assume that $\rho$ admits moments of any order: $c_n = \int_a^b t^n \rho(t)dt$ and we suppose that the space generated by the polynomials is dense in $L^2(I,\rho)$.

Stieltjes transformation of the measure of density $\rho$ is defined on the complex plane without the real interval $I$ by the formula: $z \mapsto S_\rho(z) = \int_a^b \frac{\rho(t)dt}{z-t}$. (see [1], [2]).

Let us recall the results below :

_ Let be a positive measure on $I$ associated to a density function $\mu$, also allowing moments of any order and having Stietjes's transformation linked to these of $\rho$ by the equality: $S_\mu(z) = z - c_1 - \frac{1}{S_\rho(z)}$. Secondary polynomials $Q_n$ relative with $\rho$ then form an orthogonal family for the inner product induced by $\mu$.

_ If the density ρ is continuous and provided the existence of $\varphi(x) = \lim_{\varepsilon \to 0^+} 2\int_a^b \frac{(x-t)\rho(t)dt}{(x-t)^2 + \varepsilon^2}$, we can make μ explicit by the formula : $\mu(x) = \dfrac{\rho(x)}{\dfrac{\varphi^2(x)}{4} + \pi^2\rho^2(x)}$. [see 3]

_ The operator $f(x) \mapsto g(x) = \int_a^b \dfrac{f(t) - f(x)}{t - x}\rho(t)dt$ creating secondary polynomials extends to a continuous linear map $T_\rho$ linking the space $L^2(I,\rho)$ to the Hilbert'space $L^2(I,\mu)$ whose restriction to the hyperplane $H_\rho$ of the orthogonal polynomials for ρ with $P_0 = 1$ constitutes an isometric function for both norms respectively [3]. We infer from this the formula of covariance presented below, applying to any couple (f, g) of elements of $L^2(I,\rho)$:

(1.1) $\boxed{<f/g>_\rho - <f/1>_\rho \times <g/1>_\rho = <T_\rho(f)/T_\rho(g)>_\mu}$

**Definitions.** Under the mentioned assumptions, the measure of density μ will be call *secondary measure* associated with ρ. The function φ presented above will be call the *reducer* of ρ.

**Examples.** (see [4])

- Lebesgue's measure of constant density $\rho(x) = 1$ on [0, 1]. Its reducer is defined by $\varphi(x) = 2\ln(\dfrac{x}{1-x})$.
- Chebychev's measure of the second kind $\rho(x) = \dfrac{2}{\pi}\sqrt{1-x^2}$ on [-1, 1]. Its reducer is given by $\varphi(x) = 4x$ and its secondary measure is $\mu = \dfrac{\rho}{4}$.
- The measure of density $\rho(x) = e^{-x}$ on $[0,+\infty[$, with its reducer explicated by $\varphi(x) = 2[\ln(x) - \int_0^{+\infty} e^{-t} \ln|x-t|dt] = 2e^{-x}\text{Ei}(x)$. (Ei is the exponential integral).
- Gaussian density $\rho(x) = \dfrac{e^{\frac{-x^2}{2}}}{\sqrt{2\pi}}$ on R, with $\varphi(x) = \sqrt{2\pi} e^{\frac{-x^2}{2}} \text{erfi}(\dfrac{x}{\sqrt{2}})$. ( erfi is the Imaginary Error function).

In case where ρ satisfy a Lipschitz condition on the standard interval [0, 1], we can easily explicit the reducer by: $\varphi(x) = 2\rho(x)\ln(\dfrac{x}{1-x}) - 2\int_0^1 \dfrac{\rho(t) - \rho(x)}{t-x}dt$.

If ρ is a $C^1$ class function, we have also $\varphi(x) = 2[\int_0^1 \rho'(t)\ln\left|1 - \dfrac{t}{x}\right|dt + \rho(1)\ln(\dfrac{x}{1-x})]$.

More generally, if the density $\rho$ can be evaluate at the bounds $a,b$ of the interval $I$ :

$$\varphi(x) = 2[\int_a^b \rho'(t) \ln\left|1 - \frac{t}{x}\right| dt + \rho(b) \ln(\frac{x}{b-x}) + \rho(a) \ln(\frac{x-a}{x})]$$

Other explicit examples when $I = (0, 1)$ :

- $\rho(x) = -\ln(x)$ ; $\varphi(x) = \frac{2\pi^2}{3} - \ln^2(x) - 2\text{dilog}(1-x)$

- $\rho(x) = (a+1)x^a \quad a > -1$, $a$ non integer ; $\varphi(x) = 2(a+1)[\frac{\pi x^a}{\tan(a\pi)} + \text{LerchPhi}(x,1,-a)]$

- $\rho(x) = \ln(1 + \frac{1}{\sqrt{x}})$ ; $\varphi(x) = -2\ln(\frac{1}{\sqrt{x}} + 1)\ln(\frac{1}{\sqrt{x}} - 1)$

- $\rho(x) = \frac{2}{\pi}\arccos(1-2x)$ ; $\varphi(x) = -4\ln(4-4x)$

- $\rho(x) = \frac{e^x}{e-1}$ ; $\varphi(x) = \frac{-2e^x}{e-1}[\text{Ei}(1-x) - \text{Ei}(-x)]$

Now we can present the important role of the reducer $\varphi$, in relation with $T_\rho$, and principally the inversion of this operator by an explicit formula, without using the classical decomposition in Hilbert's bases.

## 2. Notion of reducible measure and first theorem.

The measure $\rho$ will be called *reducible* if and only we have $\varphi(x) = \lim_{\varepsilon \to 0^+} 2\int_a^b \frac{(x-t)\rho(t)dt}{(x-t)^2 + \varepsilon^2}$ with the meaning of the norm of $L^2(I,\rho)$ and if the quotient $\frac{\rho}{\mu}$ belongs to the space image of $T_\rho$.

In case the space generated by the polynomials is dense in $L^2(I,\mu)$, which is the case for example when the support of the measure is compact, it is obvious that the operator $T_\rho$ is surjective. Given the relation $\frac{\rho}{\mu} = \frac{\varphi^2}{4} + \pi^2 \rho^2$, the existence of an antecedent of this quotient for $T_\rho$ still amounts to establishing that both functions $\varphi$ and $\rho$ are elements of $L^2([0,1],\rho)$.

The four first examples given above correspond to reducible measures. The Jacobi's measure of density $\rho(x) = \frac{2}{\pi}\sqrt{\frac{1-x}{x}}$ on $]0,1[$ and the Chebychev's first kind measure defined by $\rho(x) = \frac{1}{\pi\sqrt{1-x^2}}$ on $]-1, 1[$ are not reducible because $\frac{\rho}{\mu}$ is not an element of $L^2([0,1],\mu)$.

The interesting point in the notion will be specified with the following results:

**2.1. Théorem**. If the measure $\rho$ is reducible, then the unique antecedent of the quotient $\frac{\rho}{\mu}$ for the operator $T_\rho$ belonging to hyperplane $H_\rho$ is nothing but the reducer $\varphi$ defined in the section of the first theorem.

**Proof**. Let us call $\psi$ the element of $H_\rho$ sole antecedent of this quotient for $T_\rho$. We can check that $\psi = \varphi$ by establishing the coincidence of the moments of any order for the measure $\rho$, this, given the hypothesis of the density of the space of polynomial functions. Let us first evaluate the moments of $\psi$. Since $\psi$ is an element of $H_\rho$, the formula (1.1) of covariance simplifies into:

$$\sigma_n = \int_0^1 \psi(t)t^n\rho(t)dt = <\psi/x^n>_\rho = <T_\rho(\psi)/T_\rho(x^n)>_\mu = <\frac{\rho}{\mu}/T_\rho(x^n)>_\mu = <1/T_\rho(x^n)>_\rho$$

Now, an obvious calculation gives: $T_\rho(x^n) = \sum_{k=0}^{k=n-1} c_{n-1-k} x^k$, therefore

$$\sigma_n = \int_a^b \psi(t)t^n\rho(t)dt = \sum_{k=0}^{k=n-1} c_k c_{n-1-k}.$$ Let us calculate now the moments $\sigma'_n$ of $\varphi$.

We must keep in mind that $\varphi(x) = \lim_{\varepsilon \to 0^+} 2\int_a^b \frac{(x-t)\rho(t)dt}{(x-t)^2 + \varepsilon^2}$ in the space $L^2(I,\rho)$. Since the function $x \mapsto x^n$ is also an element of $L^2(I,\rho)$, we can infer with Cauchy Schwartz that

$$\sigma'_n = \lim_{\varepsilon \to 0^+} \int_a^b \left(2\int_a^b \frac{(x-t)\rho(t)dt}{(x-t)^2 + \varepsilon^2}\right) x^n\rho(x)dx.$$

Given any integer $n$, let us write for $\varepsilon$ positive and fixed:

$$\sigma'_n(\varepsilon) = 2\int_a^b\int_a^b \frac{(x-t)x^n\rho(t)dt}{(x-t)^2 + \varepsilon^2}\rho(x)dx.$$

By changing the couple $(x, t)$ into $(t, x)$ and according to Fubini we obtain:

$$\sigma'_n(\varepsilon) = \int_a^b\int_a^b \frac{(x-t)(x^n - t^n)\rho(t)\rho(x)dtdx}{(x-t)^2 + \varepsilon^2}.$$

The functions at work are now positive, hence through monotonous convergence:

$$\sigma'_n = \lim_{\varepsilon \to 0^+} \sigma'_n(\varepsilon) = \int_a^b\int_a^b \frac{x^n - t^n}{x-t}\rho(t)\rho(x)dtdx.$$ Let us now use $x^n - t^n = (x-t)\sum_{k=0}^{k=n-1} t^k x^{n-1-k}$.

By writing under the form of an integral $\sigma_n = \sum_{k=0}^{k=n-1} c_k c_{n-1-k} = \sum_{k=0}^{k=n-1} \int_a^b t^k \rho(t)dt \int_a^b x^{n-1-k}\rho(x)dx$.

We naturally come to: $\sigma_n = \int_a^b\int_a^b \frac{x^n - t^n}{x-t}\rho(t)\rho(x)dtdx = \sigma'_n$.

**Corollaries.** (Considering the theorem's hypotheses).

**2.2.** For any element $f$ of $L^2(I,\rho)$ : $<f/\varphi>_\rho = <T_\rho(f)/1>_\rho$.

By applying the formula of covariance to the couple $(f,\varphi)$ and since $<\varphi/1>_\rho = 0$ through $\varphi$ belonging to $H_\rho$, we come up with:

$$<f/\varphi>_\rho = <T_\rho(f)/T_\rho(\varphi)>_\mu = <T_\rho(f)/\frac{\rho}{\mu}>_\mu = <T_\rho(f)/1>_\rho.$$

**2.3** $\int_a^b \varphi^2(x)\rho(x)dx = \frac{4\pi^2}{3}\int_a^b \rho^3(x)dx$.

By applying the formula above with $f = \varphi$ we get, still thanks to $T_\rho(\varphi) = \frac{\rho}{\mu}$:

$$<\varphi/\varphi>_\rho = <T_\rho(\varphi)/T_\rho(\varphi)>_\mu = <\frac{\rho}{\mu}/\frac{\rho}{\mu}>_\mu = <\frac{\rho}{\mu}/1>_\rho = <\frac{\varphi^2}{4}+\pi^2\rho^2/1>_\rho$$

Hence, through linearity, in $L^2(I,\rho)$: $\|\varphi\|^2 = \frac{1}{4}\|\varphi\|^2 + \pi^2\|\rho\|^2$, and so: $\|\varphi\|^2 = \frac{4\pi^2}{3}\|\rho\|^2$

This formula easily extends by linearity to a triplet density and their associated reducer as:

$$\int_0^1 [\varphi_1(x)\varphi_2(x)\rho_3(x) + \varphi_1(x)\varphi_3(x)\rho_2(x) + \varphi_2(x)\varphi_3(x)\rho_1(x)] = 4\pi^2\int_0^1 \rho_1(x)\rho_2(x)\rho_3(x)dx$$

**2.4.** Fourier's coefficients of the reducer $\varphi$ related to an orthogonal system of polynomials $P_n$ are obtained using the associated secondary polynomials $Q_n$ by the formula:

$$C_n(\varphi) = <\varphi/P_n>_\rho = <T_\rho(\varphi)/Q_n>_\mu = <\frac{\rho}{\mu}/Q_n>_\mu = <Q_n/1>_\rho$$

**2.5.** If $f$ of class $C^1$ on $I$, with Fourier's coefficients $C_n(f)$ in an orthogonal system for $\rho$, we then have: $\int_a^b\int_a^b \frac{f(y)-f(x)}{y-x}\rho(x)\rho(y)dxdy = \sum_{n=1}^{n=\infty} C_n(\varphi)C_n(f)$

This only explicits formula 1.1. It increases the interest in the explicitation of the Fourier's coefficient of the reducer. Among the examples quoted above we get:

- For Lebesgue's measure $\rho(x)=1$ on $[0, 1]$. If we consider the Legendre polynomials conventionally defined by $P_n(x) = \frac{\sqrt{2n+1}}{n!}\frac{d^n}{dx^n}(x^n(1-x)^n)$, these coefficients are null for an even index and given by $C_n(\varphi) = -\frac{4\sqrt{2n+1}}{n(n+1)}$ for $n$ odd.

- For the Gaussian density $\rho(x) = \dfrac{e^{-\frac{x^2}{2}}}{\sqrt{2\pi}}$ on $\mathbb{R}$, if we take Hermite polynomials defined by $P_n(x) = \dfrac{1}{\sqrt{n!}} e^{\frac{x^2}{2}} \dfrac{d^n}{dx^n}(e^{-\frac{x^2}{2}})$, the $C_n(\varphi)$ are also equal to 0 for $n$ even and become explicit with $C_n(\varphi) = (-1)^{\frac{n+1}{2}} \dfrac{(\frac{n-1}{2})!}{\sqrt{n!}}$ for $n$ odd.

- For the density $\rho(x) = e^{-x}$ on $[0, +\infty[$, in relation to Laguerre's polynomials $L_n(x) = \dfrac{e^x}{n!} \dfrac{d^n}{dx^n}(x^n e^{-x}) = \sum_{k=0}^{k=n} C_n^k (-1)^k \dfrac{x^k}{k!}$, the Fourier's coefficients of the reducer are given by : $C_n(\varphi) = -\dfrac{1}{n} \sum_{k=0}^{k=n-1} \dfrac{1}{C_{n-1}^k}$.

This last example deserves further development. First of all we notice that the calculation of these coefficients naturally shows the summation of the elements of the line indexed by $n$ in the table of the Leibniz's harmonic triangular numbers. A simple application of the formula 2.3 informs us on the summation of these constants squared.

$$\int_0^{+\infty} 4[\text{Ei}(x)]^2 e^{-3x} dx = \dfrac{4\pi^2}{9} = \sum_{n=1}^{n=+\infty} \left( \dfrac{1}{n} \sum_{k=0}^{k=n-1} \dfrac{1}{C_{n-1}^k} \right)^2$$

Puzzled by this unexpected triangle, I went further into the study of these summations $s_n = \dfrac{1}{n} \sum_{k=0}^{k=n-1} \dfrac{1}{C_{n-1}^k}$. (see [5])

Indeed we come up with a far more pleasant expression noted : $s_n = \dfrac{1}{2^n} \sum_{k=1}^{k=n} \dfrac{2^k}{k}$.

I have no knowledge of a standard designation for numbers $G_n = \sum_{k=1}^{k=n} \dfrac{2^k}{k}$.

It seems to me that naming them 'geoharmonic' would be appropriate referring to the triangle in question of course and the powers of 2 but also because of a duality in formulas involving these values and the conventional harmonic numbers (see [6] ). The table below speaks for itself.

| **Harmonic numbers** $H_n = \sum_{k=1}^{k=n} \frac{1}{k}$ | **'Geoharmonic' numbers** $G_n = \sum_{k=1}^{k=n} \frac{2^k}{k}$ |
|---|---|
| $\sum_{n=1}^{n=\infty} H_n x^n = \frac{\ln(1-x)}{x-1}$ | $\sum_{n=1}^{n=\infty} G_n x^n = \frac{\ln(1-2x)}{x-1}$ |
| $\sum_{n=1}^{n=\infty} \frac{H_n}{n} x^n = \frac{1}{2}\ln^2(1-x) + \text{dilog}(1-x)$ | $\sum_{n=1}^{n=\infty} \frac{G_n}{n} x^n = \frac{1}{2}\ln^2(1-x) + \text{dilog}(1-x) + \text{dilog}(\frac{1-2x}{1-x})$ |
| $\sum_{n=1}^{n=\infty} \frac{H_n}{n!} x^n = e^x \int_0^x \frac{1-e^{-t}}{t} dt = e^x \sum_{n=1}^{n=\infty} \frac{(-1)^{n+1} x^n}{n \cdot n!}$ | $\sum_{n=1}^{n=\infty} \frac{G_n}{n!} x^n = e^x \int_0^x \frac{e^t - e^{-t}}{t} dt = e^x \sum_{n=1}^{n=\infty} \frac{[1+(-1)^{n+1}] x^n}{n \cdot n!}$ |
| $\sum_{n=1}^{n=\infty} \frac{H_n^2}{4^n} = \frac{2\pi^2}{9} + \frac{8}{3}[\ln(2)\ln(\frac{2}{3}) - \text{dilog}(\frac{2}{3})]$ | $\sum_{n=1}^{n=\infty} \frac{G_n^2}{4^n} = \frac{4\pi^2}{9}$ |

I also obtain the equality: $\boxed{\sum_{n=1}^{n=\infty} \frac{G_n H_n}{4^n} = \frac{4}{3}[\frac{\pi^2}{6} + \text{dilog}(\frac{3}{2}) - \ln(2)\ln(\frac{3}{2})]}$.

(For the dilogarithme function see [7], [8]).

What is even more interesting is the link appearing between the continuous density $\rho(x) = e^{-x}$ and the discrete measure of weight $p_n = \frac{1}{2^n}$.

By studying the Fourier's coefficients of the product of the reducer $\varphi$ by a monomonial we show the equality: $C_n^k(\varphi) = \int_0^{+\infty} \varphi(x) x^k L_n(x) e^{-x} dx = A_k(n) s_n - B_k(n)$, with $k \mapsto A_k$ a sequence of orthogonal polynomials for the discrete weight $p_n = \frac{1}{2^n}$ and $B_k$ secondary polynomial associated with $A_k$ for that same weight. The demonstration is carried out by analysing the relation of recurrence.

## 3. Inverted operator.

Let us suppose that the support of the measure is a compact interval and to make things easier we will work on the standard interval $I = [0,1]$.

**3.1. Theorem.** Let us assume the measure $\rho$ is reducible. For any element $f$ of $L^2([0,1],\rho)$ that $\varphi \times f$ et $T_\rho(f)$ belong to $L^2([0,1],\rho)$, and $\frac{\rho}{\mu} \times f$ is an element of $L^2([0,1],\mu)$, we get the so-called composition formula: $\boxed{T_\rho(\varphi \times f - T_\rho(f)) = \frac{\rho}{\mu} \times f}$

**Proof.**

Under the preceding hypotheses, the function $h = \varphi \times f - T_\rho(f)$ is an element of $L^2(I,\rho)$, and verifies $\int_0^1 h(t)\rho(t)dt = 0$. Indeed according to (1.1) we can write:

$<\varphi \times f / 1>_\rho = <\varphi / f>_\rho = <T_\rho(f)/1>_\rho$. The function $h$ is therefore an element of the hyperplane $H_\rho$ and consequently $T_\rho(h)$ belongs to $L^2(I,\mu)$.

The difference $T_\rho(h) - \dfrac{\rho}{\mu} \times f$ is then an element of $L^2([0,1],\mu)$.

To show its null value we only have to proof it is orthogonal to any element in an orthogonal basis of the Hilbert's space $L^2([0,1],\mu)$.

If we keep $n \mapsto P_n$ to refer to the initial system of orthogonal polynomials relating to $\rho$, we just have to state that for any integer $n$ :

$<T_\rho(\varphi \times f - T_\rho(f))/T_\rho(P_n)>_\mu = <\dfrac{\rho}{\mu} \times f / T_\rho(P_n)>_\mu$.

The first inner product is worth the formula of covariance, and as $h$ belongs to $H_\rho$ :

$<\varphi \times f - T_\rho(f)/P_n>_\rho = <\varphi/f \times P_n>_\rho - <T_\rho(f)/P_n>_\rho$.

We can also write : $<\varphi/f \times P_n>_\rho = <\dfrac{\rho}{\mu}/T_\rho(f \times P_n)>_\mu = <T_\rho(f \times P_n)/1>_\rho$

The second scalar product simplifies into $<f \times T_\rho(P_n)/1>_\rho$. What we have to do now is to check :

$$\boxed{\int_0^1 T_\rho[f(t)P_n(t)]\rho(t)dt = \int_0^1 T_\rho(f(t))P_n(t)\rho(t)dt + \int_0^1 f(t)T_\rho(P_n(t))\rho(t)dt} \quad (F)$$

We will first establish it when $f$ is a polynomial function.

Let us first check it in the particular case of $f(x) = x$.

It is easy to show : $T_\rho(tP_n(t)) = tT_\rho(P_n(t)) + \int_0^1 P_n(t)\rho(t)dt$

We now have to multiply by $\rho(t)$ and integrate on [0,1], which results in the expected equality (F) in this particular case for $\int_0^1 \rho(t)dt = c_0 = 1$ and $T_\rho(x) = 1$.

Let us now consider any whole power $f(x) = x^k$.

It is easy to verify by recurrence on $k$ that for all polynomial $P$ :

$$T_\rho(x^k P(x)) = x^k T_\rho(P(x)) + \sum_{j=0}^{j=k-1} C(P,k,j)x^j \quad \text{with } C(P,k,j) = \int_0^1 t^{k-1-j} P(t)\rho(t)dt$$

Particularly : $T_\rho(x^k) = \sum_{j=0}^{j=k-1} c_{k-1-j} x^j$.

Therefore infer that : $\int_0^1 T_\rho(t^k P_n(t))\rho(t)dt = \int_0^1 t^k T_\rho(P_n(t))\rho(t)dt + \sum_{j=0}^{j=k-1} C(P_n,k,j) \times c_j$.

Now, given the formula above of the transformation of a power :

$$\int_0^1 T_\rho(t^k) P_n(t)\rho(t)dt = \sum_{j=0}^{j=k-1} c_{k-1-j} \int_0^1 t^j P_n(t)\rho(t)dt = \sum_{j=0}^{j=k-1} c_{k-1-j} \times C(P_n,k,k-1-j) = \sum_{j=0}^{j=k-1} c_j \times C(P_n,k,j)$$

We still obtain the equality (*F*) for any whole power function and subsequently by linearity of the integral, for any polynomial function *f*.

Finally, to establish the validity of (*F*) for any element *f* that satisfies the restrictive conditions mentioned, we only have to consider it as the limit of a sequence of polynomials $f_k$ while the meaning of the norm on L²([0,1], ρ )

As the function $\frac{\rho}{\mu}$ is an element of L²([0,1], μ ), we will have no difficulty in passing to the limit for the three integrals considered.

Let us analyse in full details these delicate operations.

➢ For the left hand member in (*F*). If we write $g_k(t) = f(t) - f_k(t)$, we will have:

$$I_k = \int_0^1 T_\rho[f(t)P_n(t)]\rho(t)dt - \int_0^1 T_\rho[f_k(t)P_n(t)]\rho(t)dt = \int_0^1 T_\rho[g_k(t)P_n(t)]\rho(t)dt$$

As assumed, the difference $g_k$ tends to 0 in the space L²([0,1],ρ). Polynomial $P_n$ being bounded on [0,1], the product $g_k \times P_n$ also tends to 0 for the same norm, and its image by the continuous map $T_\rho$ will tend to 0 for the norm on L²([0,1],μ )

Let us write $I_k = \int_0^1 T_\rho[g_k(t)P_n(t)]\frac{\rho(t)}{\mu(t)}.\mu(t)dt$ ; Since $\frac{\rho}{\mu}$ is an element of L²([0,1],μ ), $I_k$ turns out to be a scalar product in this space and thus can be majored in absolute value with Cauchy-Schwartz. $|I_k| \leq \|T_\rho(g_k P_n)\|_\mu \times \left\|\frac{\rho}{\mu}\right\|_\mu$ We then deduce $\lim_{k \to \infty} I_k = 0$.

➢ As concerns the first integral of the right-hand member (*F*) we will follows the same steps.

$$J_k = \int_0^1 T_\rho[f(t)]P_n(t)\rho(t)dt - \int_0^1 T_\rho[f_k(t)]P_n(t)\rho(t)dt = \int_0^1 T_\rho[g_k(t)]P_n(t)\rho(t)dt$$

$|J_k| \leq \|T_\rho(g_k)\|_\mu \times \|P_n \times \dfrac{\rho}{\mu}\|_\mu$. We also conclude with $\lim\limits_{k \to \infty} J_k = 0$.

➢ Finally for the last integral in (*F*), the conclusion is more simple :

$$L_k = \int_0^1 f(t) T_\rho(P_n(t)) \rho(t) dt - \int_0^1 f_k(t) T_\rho(P_n(t)) \rho(t) dt = \int_0^1 g_k(t) T_\rho(P_n(t)) \rho(t) dt$$

$$|L_k| \leq \|g_k\|_\rho \times \|T_\rho(P_n)\|_\rho \ .$$

Indeed $T_\rho(P_n) = Q_n$ is a polynomial, therefore an element of $L^2([0,1], \rho)$. So, we come to the same $\lim\limits_{k \to \infty} L_k = 0$.

Subsequently the theorem is clearly established. If the function *f* of $L^2([0,1], \rho)$ is such that $\varphi \times f$ and $T_\rho(f)$ also are elements of $L^2([0,1], \rho)$ and such $\dfrac{\rho}{\mu} \times f$ be an element of $L^2([0,1], \mu)$ we can write :

$$\boxed{T_\rho(\varphi \times f) = (T_\rho)^2(f) + \dfrac{\rho}{\mu} \times f}$$

**Please note.**

The assumptions of the theorem are very constraining in their generality but very easy to verify for instance in the case of a function of class $C^1$ on [0, 1], considering bounded measures.

**Corollaries.**

**3.2.** If the measure $\rho$ is reducible, we can define an isometry $U_\rho$ from the space $L^2([0,1], \dfrac{\rho^2}{\mu})$ towards the hyperplane $H_\rho$ provided with the norm of $L^2([0,1], \rho)$, satisfying the relation of composition : $T_\rho \circ U_\rho(f) = \dfrac{\rho}{\mu} \times f$.

**Proof.** According to the stated theorem and if *f* is a polynomial, we can write with $h = U(f) = \varphi \times f - T_\rho(f)$, that $T_\rho(h) = \dfrac{\rho}{\mu} \times f$.

Now, *h* being an element of $H_\rho$, we infer from the formula of covariance (1.1):

$$<h/h>_\rho = <T_\rho(h)/T_\rho(h)>_\mu = <\dfrac{\rho}{\mu} \times f / \dfrac{\rho}{\mu} \times f>_\mu = <f/f>_{\frac{\rho^2}{\mu}}$$

We then obtain the identity: $\boxed{\int_0^1 h^2(t) \rho(t) dt = \int_0^1 f^2(t) \dfrac{\rho^2(t) dt}{\mu(t)}}$

This equality can be written in a simple way : $\boxed{\|U(f)\|_\rho = \|f\|_{\frac{\rho^2}{\mu}}}$

So it will be possible to extend the operator $U$, according to Cauchy's theorem, into an isometry $U_\rho$ from the space $L^2([0,1], \frac{\rho^2}{\mu})$ towards $H_\rho$ provided with the norm of $L^2([0,1], \rho)$. (We shall not forget here that $U_\rho$ is not surjective a-priori because the images of the polynomials are not polynomials of an orthogonal basis of the arrival space).

However we will have from any element $f$ of $L^2([0,1], \frac{\rho^2}{\mu})$ : $T_\rho \circ U_\rho(f) = \frac{\rho}{\mu} \times f$

Indeed, if $f$ is a limit of a sequence of polynomials $f_k$ in the meaning of the norm of $L^2([0,1], \frac{\rho^2}{\mu})$, we can infer that $U_\rho(f_k)$ tends to $U_\rho(f)$ in the meaning of the norm of $L^2([0,1], \rho)$ and consequently $T_\rho(U_\rho(f_k))$ tends to $T_\rho(U_\rho(f))$ in the meaning of the norm of $L^2([0,1], \mu)$.

Now we know that $T_\rho(U_\rho(f_k)) = \frac{\rho}{\mu} \times f_k$, for $f_k$ is a polynomial, and it tends to $\frac{\rho}{\mu} \times f$ in $L^2([0,1], \mu)$. This simply because : $\int_0^1 [\frac{\rho}{\mu}(f - f_k)]^2 \mu(x)dx = \int_0^1 (f - f_k)^2 \frac{\rho^2(x)}{\mu(x)} dx$

We actually have the announced equality.

To come to an end, we will check that when $f$ is an element of $L^2([0,1], \frac{\rho^2}{\mu})$ meeting the restrictive conditions of the theorem of composition (3.1) above, then $U_\rho$ interferes on $f$ like : $U_\rho(f) = U(f) = \varphi \times f - T_\rho(f)$.

This is quite clear since given the theorem, in this case : $T_\rho(U(f)) = \frac{\rho}{\mu} \times f$.

So we get : $T_\rho(U(f)) = \frac{\rho}{\mu} \times f = T_\rho(U_\rho(f))$. We conclude with the injectivity of $T_\rho$ on $H_\rho$.

Let us finally notice that $U_\rho$ acts as an adjunct of $T_\rho$ for the inner product induced by $\rho$ for, on conditions that the considered compositions are possible, we can write :

$$\boxed{<U_\rho(f)/g>_\rho = <T_\rho(U_\rho(f))/T_\rho(g)>_\mu = <\frac{\rho}{\mu} \times f / T_\rho(g)>_\mu = <f/T_\rho(g)>_\rho}$$

**3.3.** Under the hypotheses and with the above notations, the operator $U_\rho$ is bijective, just as the restriction $\tilde{T}_\rho$ of $T_\rho$ to the hyperplane $H_\rho$ and we explicit their reciprocal elements thanks to the formula of composition introduced above.

**Proof**. Let us write $F_{(\rho,\mu)}$ the function from $L^2([0,1], \frac{\rho^2}{\mu})$ towards $L^2([0,1],\mu)$ changing $f$ into $\frac{\rho}{\mu} \times f$. We obviously have here an isometry an its reciprocal element is nothing that the function we will write $F_{(\mu,\rho)}$ defined by : $f \mapsto F_{(\mu,\rho)}(f) = \frac{\mu}{\rho} \times f$. So we can translate the formula of composition by a coherent diagram linking $\tilde{T}_\rho$ and $U_\rho$ : $\boxed{\tilde{T}_\rho \circ U_\rho = F_{(\rho,\mu)}}$.

The application $\tilde{T}_\rho$ restricted to $H_\rho$ being bijective, like $F_{(\rho,\mu)}$, we can infer that $U_\rho$ is bijective too. This gives us: $\boxed{(U_\rho)^{-1} = F_{(\mu,\rho)} \circ \tilde{T}_\rho}$ and $\boxed{(\tilde{T}_\rho)^{-1} = U_\rho \circ F_{(\mu,\rho)}}$.

These relations are very convenient to solve the integral equation [see 9] of type :

$$\int_0^1 \frac{f(t) - f(x)}{t - x} \rho(t) dt = g(x)$$, with $g$ being continuous on [0, 1] and $f$ unknown in $H_\rho$.

What is stated above shows that the function defined by :

$$f(x) = \varphi(x) \times \frac{\mu(x)}{\rho(x)} \times g(x) - T_\rho(\frac{\mu(x)}{\rho(x)} \times g(x))$$ is a solution on condition that $\frac{\mu}{\rho} \times g$ answers the hypothesis of (3.1), which is the case for example if this function is of class $C^1$.

**3.4.** In fact it can be shown that the expression of the solution above is simplified as:

$$f(x) = (x - c_1) g(x) - T_\mu(g(x))$$, with $c_1$ moment of order 1 for the density $\rho$.

We give in what follows a simplified proof of this when the support of the measure $\rho$ is a compact $I = [0,1]$.

We first study the reducer of the secondary measure of density $\mu$.

Recall the formula of coupling : $S_\mu(z) = z - c_1 - \dfrac{1}{S_\rho(z)}$

It is known from Stieltjes –Perron : $\mu(x) = \lim_{\varepsilon \to 0^+} \dfrac{S_\mu(x - i\varepsilon) - S_\mu(x + i\varepsilon)}{2i\pi}$

The reducer of $\mu$ is defined by: $\psi(x) = \lim_{\varepsilon \to 0^+} S_\mu(x - i\varepsilon) + S_\mu(x + i\varepsilon)$

This gives us : $\psi(x) = 2(x - c_1) - \lim_{\varepsilon \to 0^+} \dfrac{S_\rho(x - i\varepsilon) + S_\rho(x + i\varepsilon)}{S_\rho(x - i\varepsilon) \times S_\rho(x + i\varepsilon)}$

So : $\psi(x) = 2(x - c_1) - \dfrac{\varphi(x)}{\dfrac{\varphi^2(x)}{4} + \pi^2 \rho^2(x)}$ and simply we obtain:

(3.4.1) $\psi(x) = 2(x - c_1) - \dfrac{\varphi(x)\mu(x)}{\rho(x)}$

Recall now that in the case where ρ satisfy a Lipschitz condition on the standard interval $I = [0,1]$, we can easily explicit the reducer by: $\varphi(x) = 2\rho(x)\ln(\frac{x}{1-x}) - 2\int_0^1 \frac{\rho(t) - \rho(x)}{t - x} dt$.

So we can write in simplified notations under above conditions:

$$\varphi(x) = 2\rho(x)\ln(\frac{x}{1-x}) - 2T_1(\rho)$$

$$\psi(x) = 2\mu(x)\ln(\frac{x}{1-x}) - 2T_1(\mu)$$

And consequently by simple combinations

(3.4.2) $\quad \varphi(x)\mu(x) - \psi(x)\rho(x) = 2[\rho(x)T_1(\mu) - \mu(x)T_1(\rho)]$

We can now prove the identity

(3.4.3) $\quad \varphi(x) \times \frac{\mu(x)}{\rho(x)} \times g(x) - T_\rho(\frac{\mu(x)}{\rho(x)} \times g(x)) = (x - c_1)g(x) - T_\mu(g(x))$

An equivalent form with simplified notations is $T_\rho(\frac{\mu.g}{\rho}) - T_\mu(g) = g(x).[\frac{\mu(x).\varphi(x)}{\rho(x)} - (x - c_1)]$

It is very easy to shown that : $T_\rho(\frac{\mu.g}{\rho}) - T_\mu(g) = \frac{g(x)}{\rho(x)} \int_0^1 \frac{\mu(t)\rho(x) - \mu(x)\rho(t)}{t - x} dt$

Note then : $\int_0^1 \frac{\mu(t)\rho(x) - \mu(x)\rho(t)}{t - x} dt = \rho(x)T_1(\mu) - \mu(x)T_1(\rho)$

Thus, the identity (3.4.3) is equivalent to : $\frac{\mu(x)\varphi(x)}{\rho(x)} - (x - c_1) = T_1(\mu) - \frac{\mu(x)T_1(\rho)}{\rho(x)}$, and also equivalent thanks (3.4.2) to : $2[\mu(x)\varphi(x) - (x - c_1)\rho(x)] = \varphi(x)\mu(x) - \psi(x)\rho(x)]$ , who is a direct consequence of (3.4.1). This concludes the proof.

**References.**

[1] Christian Berg, Annales de la Faculté des Sciences de Toulouse, Sér 6 Vol. S5 (1996), p.9-32.

[2] G.A Baker, Jr and P. Graves-Morris, 'Padé approximants', Cambridge University Press, London 1996.

[3] R. Groux. C.R Acad.Sci. Paris, Ser.I. 2007. Vol 345. pages 373-376.

[4] A.Nikiforov, V. Ouvarov, Elements de la théorie des fonctions spéciales, Editions de Moscou, 1983.

[5]  Sury, European Journal of Combinatorics, 1993, 14, 351-353

[6] J.H Conway and R.K Guy, The book of Numbers, New York, Springer-Verlag, 1996.

[7] M Abramowitz and I.A Stegun, Handbook of Mathematical Functions, National Bureau of Standards, 1964, reprinted Dover Publications, 1965.

[8] L. Lewin, Dilogarithms and Associated Functions, London, Macdonald, 1958.

[9] A.D Polyanin and A.V Manzhirov, Handbook of integral equations, CRC Press, Boca Raton, 1998.